# EXISTENTIALLY CLOSED MODELS OF THE
# THEORY OF ARTINIAN LOCAL RINGS

Hans Schoutens

21.05.97

ABSTRACT. The class of all Artinian local rings of length at most $l$ is $\forall_2$-elementary, axiomatised by a finite set of axioms $\mathcal{A}\mathfrak{rt}_l$. We show that its existentially closed models are Gorenstein, of length exactly $l$ and their residue fields are algebraically closed, and, conversely, every existentially closed model is of this form. The theory $\mathcal{G}\mathfrak{or}_l$ of all Artinian local Gorenstein rings of length $l$ with algebraically closed residue field is model complete and the theory $\mathcal{A}\mathfrak{rt}_l$ is companionable, with model-companion $\mathcal{G}\mathfrak{or}_l$.

## 0. INTRODUCTION

In studying classes of commutative rings with identity (*rings*, for short) from a model-theoretic point of view, the elementary or (first order) axiomatisable classes are the most interesting. Given a general property of rings, one might look for elementary classes all of whose models share this property. Among their beloved properties, commutative algebraists will no doubt quote the Noetherian property. Unfortunately, this is (most of the time) a non-first order property: one can prove using ultraproducts, that any elementary class containing at least one ring with a non-zero divisor which is not a unit, must contain non-Noetherian rings.

Therefore, if one wants to maintain Noetherianness, one is then lead to look for elementary classes consisting entirely out of Artinian rings. However, if the length is unbounded, then again such classes cannot exist by an ultraproduct argument. One might bound the length by bounding some other invariants, for instance the embedding dimension and the exponent, but in this paper we will just bound the length itself. Artinian rings are in general semi-local, but for sake of simplicity, we will restrict us here to local rings. A simple but basic observation in this setting is that any extension of Artinian local rings is necessarily local, since the non-units are precisely the nilpotent elements.

There is a second motivation for the present research stemming from algebraic geometry. It is one of the many insights of GROTHENDIECK that even if one were merely interested in the study of irreducible varieties over an algebraically closed field, one nevertheless has to include in ones study the behaviour of nilpotent functions. This standpoint forces us to shift from classical algebraic geometry (corresponding to radical ideal theory) to the general sheaf-theory on schemes (corresponding to the study of arbitrary ideals). Therefore, if model theory seeks to enter in a serious dialog with algebraic geometry, it will have to enlarge its study of sets definable over an algebraic closed field to the study of sets definable over an Artinian local ring. Using regular sequences (a first order definable concept) then will enable us to reduce the study of arbitrary Cohen-Macaulay local rings to the study of Artinian local rings, thus providing a framework to mimic the infinitesimal study of scheme theory, at least under the Cohen-Macaulay assumption. The latter reduction will appear in a future paper in which the higher dimensional versions of the present theory will be given.

In short, we propose in this paper the study of Artinian local rings of length at most $l$ from a model-theoretic point of view, for some $l \in \mathbb{N}$. If we moreover require that the rings are equicharacteristic, then by COHEN's structure theorem, each such ring is a finite dimensional algebra over some (coefficient) field (isomorphic to the residue field). However, this field need not to be definable.[1]

It turns out that the class of Artinian local rings of length at most $l$ is axiomatisable by a $\forall_2$-theory, denoted $\mathcal{A}\mathfrak{rt}_l$. To model-theorists, $\forall_2$-theories behave particularly nicely and one of their important features is that any model is embeddable in an existentially closed model. A model $R$ is called *existentially closed* or *generic*, whenever any $\exists_1$-sentence with parameters from $R$ which holds in an extension $S$ of $R$ (i.e., such that $R \subset S$ and $S$ is also a model), then it holds already in $R$. In other words, solvability of a system of equations in $R$ can be tested in any extension of $R$. For example, the existentially closed models of the theory of fields are precisely the algebraically closed fields.

---

[1]In the case of positive characteristic $p$ and perfect residue field the coefficient field is definable as the collection of all elements in the ring which are an $p^l$-th power.

Our main theorem states that an existentially closed model of $\mathcal{A}\mathfrak{rt}_l$ is Gorenstein, has an algebraically closed residue field and its length is exactly equal to $l$ and there are no other existentially closed models. In commutative algebra or algebraic geometry, the Gorenstein condition is the 'next best property' after being regular or being a complete intersection. Gorenstein rings were initially studied by GROTHENDIECK and SERRE in algebraic geometry and by BASS in commutative algebra, see for instance the latter's paper [**Bass**]. For our purposes, we only need to define the Gorenstein property for Artinian local rings, in which case it is equivalent with being self-injective; or with having a unique non-zero ideal; or with having a principal socle (=annihilator of the maximal ideal). Therefore, Artinian local Gorenstein rings are better known to logicians as quasi-Frobenius rings. One of its less used characterisations will be needed in this paper: they are precisely those Artinian local rings $R$ for which each ideal is an *annihilator* ideal, meaning, an ideal of the form $\operatorname{Ann}_R(I)$, where $I$ is an ideal of $R$. It is precisely this characterisation together with exercise [**Hod 1**, 3.2. Exercise 17], about a similar property of existentially closed rings, that lead us to the present research.

An algebraic consequence of this theorem is that any Artinian local ring can be embedded in an Artinian local Gorenstein ring of the same length. We give an example **(1.8)** to show that this is not at all that obvious. Model-theoretic consequences are that the theory $\mathcal{G}\mathfrak{or}_l$ of all Artinian local Gorenstein rings of length $l$ with algebraically closed residue field is model-complete. It is the model-companion of $\mathcal{A}\mathfrak{rt}_l$. However, in the last section we give a counterexample to Quantifier Elimination for these theories. In a subsequent paper [**Sch 2**] we will show how one can obtain Quantifier Elimination by adding some extra constants and some natural predicates.

In order to prove that existentially closed models of $\mathcal{A}\mathfrak{rt}_l$ are models of $\mathcal{G}\mathfrak{or}_l$ (see **(1.6)**), the only place at which I really seem to need some model-theory is in showing that the length of an existentially closed model is the maximal possible value $l$. To prove the converse statement, we need some better understanding of the ideal structure of an Artinian local Gorenstein ring. In section 2 we deal with the easier case when the residue field is embeddable[2], the so called equicharacteristic case. We show the existence of certain monomial bases having the right combinatorial properties due to the Gorenstein condition. We show that equicharacteristic Artinian local Gorenstein rings are in some sense rigid, in that any extension by an Artinian local ring of the same length has to be automatically Gorenstein and, even more, is obtained by base change over the residue fields (see **(2.3)**). In particular, any such extension is automatically faithfully $\Omega$at. This result allows us to translate any system of equations over the model to a system of equations over the residue field, which is algebraically(=existentially) closed.

The same ideas are then exploited in the next section to treat the remaining mixed characteristic case. The lack of a coefficient field is now the main obstruction. However, we show how the multiplicative part of the residue field is still embeddable, using Witt vectors. Therefore a similar rigidity for extensions of a model of $\mathcal{G}\mathfrak{or}_l$ as in the equicharacteristic case still exists, but whereas the former is easily expressed by a base change, the latter requires a more careful analysis of the terms involved (see **(3.7)**) to finish the argument.

### 0.0. Definitions and Notations.

#### Logic

The standard work of reference will be the excellent book [**Hod 2**]. Most of our terminology and notations will be consistent with it. However, the numbering within tuples sometimes starts with one instead of zero. The length of a tuple $a$ will be denoted $\operatorname{len}(a)$.

Let $\mathcal{L}$ be a (first order) language. We will say that a class $\mathbb{H}$ of $\mathcal{L}$-structures is *elementary*, if there exists a first order theory $T$, such that its models form the class $\mathbb{H}$. We will say that the class $\mathbb{H}$ is $\forall_n$-*elementary*, where $n$ is some natural number, if $T$ consists entirely out of $\forall_n$-sentences.

---

[2]By COHEN's Structure Theorem we are then in the case of a finite dimensional local algebra over a field (called the *coefficient field*)

In the rest of this paper, $\mathcal{L}$ will always stand for the language of rings

$$\mathcal{L} =< +, \cdot, -; \mathbf{0}, \mathbf{1} >$$

where $+$ and $\cdot$ are two-ary function symbols interpreted by respectively addition and multiplication, where $-$ is a unary function symbol interpreted by the additive inverse, and where $\mathbf{0}$ and $\mathbf{1}$ are constant symbols, standing for the zero and one of a ring. The theory $\mathcal{R}ing$ of rings (commutative, with identity) is $\forall_1$ so that the class of all rings is $\forall_1$-elementary.

A notational convention: symbols of the language and logical variables will appear in *boldface*, whereas their equivalent in normal face will stand for their interpretation in a particular model, sometimes emphasised by adding a superscript denoting the structure. Moreover, structures will always be denoted by the same symbol as their domain, confusing as this might be.

## Commutative Algebra

Here, the work of reference will be [**Mats**] . Let $R$ be a local ring. Its residue field will normally be denoted by $\kappa$ and its maximal ideal by $\mathfrak{m}$. We might want to indicate this by writing $(R, \mathfrak{m})$ or $(R, \mathfrak{m}, \kappa)$. We will denote the *length* of a ring $R$, i.e., the length of a maximal chain of ideals in $R$, by $\ell(R)$. By $\exp(R)$ we will mean the *exponent* or *nilpotency degree* of $R$, defined as the minimum of all $e$, such that $\mathfrak{m}^e = 0$, if such an $e$ exists, otherwise we set $\exp(R) = \infty$. The *embedding dimension* of $R$ will be the cardinality of a minimal set of generators of the maximal ideal $\mathfrak{m}$. By NAKAYAMA's lemma, this equals the dimension of $\mathfrak{m}/\mathfrak{m}^2$ over the residue field $\kappa$.

**0.1. Proposition.** *The class of all local rings is $\forall_2$-elementary.*

*Proof.* We express that an element is a non-unit, by

(0.1.1) $$\mathtt{Nu}(\boldsymbol{r}) = (\forall \boldsymbol{s})[\boldsymbol{s}\boldsymbol{r} \neq \mathbf{1}].$$

Hence, the following sentence

(0.1.2) $$\mathtt{Loc} \stackrel{\text{def}}{=} (\forall \boldsymbol{r}, \boldsymbol{s})[\mathtt{Nu}(\boldsymbol{r}) \wedge \mathtt{Nu}(\boldsymbol{s}) \to \mathtt{Nu}(\boldsymbol{r} + \boldsymbol{s})]$$

expresses that a ring $R$ is local. Writing this out with aid of (0.1.1), $\mathtt{Loc}$ becomes a $\forall_2$-sentence. ∎

**0.2. Reduction Lemma.** *Let $R$ be a ring and let $\bar{\boldsymbol{a}}$ be an $n$-tuple of free variables. There exists a transformation which turns a formula $\mathtt{P}(\bar{\boldsymbol{x}})$ of $\mathcal{L}$ with parameters from $R$ in the $m$-tuple $\bar{\boldsymbol{x}}$ of free variables, into a formula $\mathrm{Red}_{\bar{\boldsymbol{a}}} \mathtt{P}(\bar{\boldsymbol{x}})$ of $\mathcal{L}$ with (the same) parameters from $R$ in the free variables $\bar{\boldsymbol{a}}$ and $\bar{\boldsymbol{x}}$, such that, for any ideal $\mathfrak{a}$ of $R$ generated by $n$ elements $a_1, \ldots, a_n$ and any $x \in R^m$, we have that*

$$R \models \mathrm{Red}_a \mathtt{P}(x) \iff R/\mathfrak{a} \models \mathtt{P}(x),$$

*where we have reduced the parameters modulo $\mathfrak{a}$ in the latter formula. Moreover, if $\mathtt{P}$ is either an $\exists_k$ formula, for some odd $k$, or an $\forall_k$ formula, for some even $k > 0$, then likewise is $\mathrm{Red}_{\bar{\boldsymbol{a}}} \mathtt{P}$.*

*Proof.* By definition, any formula is built up from the following four rules.

(F1) For each pair of terms $t_1$ and $t_2$, the expression $t_1 = t_2$ is a formula.
(F2) If $\varphi$ is a formula, then also $\neg \varphi$.
(F3) If $\varphi$ and $\psi$ are formulae, then also $\varphi \wedge \psi$.
(F4) If $\varphi$ is a formula, then also $(\exists \boldsymbol{z})\varphi$.

We will now transform $\mathtt{P}(\bar{\boldsymbol{x}})$ into a formula $\mathrm{Red}_{\bar{\boldsymbol{a}}} \mathtt{P}(\bar{\boldsymbol{x}})$ as follows, by induction on the number of formation rules (F1)–(F4) needed to construct $\mathtt{P}(\bar{\boldsymbol{x}})$.

(R1) If $t_1(\bar{\boldsymbol{x}})$ and $t_2(\bar{\boldsymbol{x}})$ are terms, and $\mathtt{P}(\bar{\boldsymbol{x}})$ stands for the formula $t_1(\bar{\boldsymbol{x}}) = t_2(\bar{\boldsymbol{x}})$, then $\mathrm{Red}_{\bar{\boldsymbol{a}}} \mathtt{P}(\bar{\boldsymbol{x}})$ stands for the formula $(\exists \bar{\boldsymbol{y}})[t_1(\bar{\boldsymbol{x}}) - t_2(\bar{\boldsymbol{x}}) = \sum_{i=1}^{n} \boldsymbol{a_i} \boldsymbol{y_i}]$.
(R2) If $\varphi(\bar{\boldsymbol{x}})$ is a formula and $\mathtt{P}(\bar{\boldsymbol{x}})$ stands for the formula $\neg \varphi(\bar{\boldsymbol{x}})$, then $\mathrm{Red}_{\bar{\boldsymbol{a}}} \mathtt{P}(\bar{\boldsymbol{x}})$ stands for the formula $\neg \mathrm{Red}_{\bar{\boldsymbol{a}}} \varphi(\bar{\boldsymbol{x}})$.
(R3) If $\varphi(\bar{\boldsymbol{x}})$ and $\psi(\bar{\boldsymbol{x}})$ are formulae and $\mathtt{P}(\bar{\boldsymbol{x}})$ stands for the formula $\varphi(\bar{\boldsymbol{x}}) \wedge \psi(\bar{\boldsymbol{x}})$, then $\mathrm{Red}_{\bar{\boldsymbol{a}}} \mathtt{P}(\bar{\boldsymbol{x}})$ stands for the formula $\mathrm{Red}_{\bar{\boldsymbol{a}}} \varphi(\bar{\boldsymbol{x}}) \wedge \mathrm{Red}_{\bar{\boldsymbol{a}}} \psi(\bar{\boldsymbol{x}})$.
(R4) If $\varphi(\bar{\boldsymbol{x}})$ is a formula and $\mathtt{P}(\bar{\boldsymbol{x}})$ stands for the formula $(\exists \boldsymbol{z})\varphi(\bar{\boldsymbol{x}}, \boldsymbol{z})$, then $\mathrm{Red}_{\bar{\boldsymbol{a}}} \mathtt{P}(\bar{\boldsymbol{x}})$ stands for the formula $(\exists \boldsymbol{z})\mathrm{Red}_{\bar{\boldsymbol{a}}} \varphi(\bar{\boldsymbol{x}}, \boldsymbol{z})$.

It is now not hard to verify that $\text{Red}_a\,\mathtt{P}(x)$ holds in $R$, if and only if, $\mathtt{P}(x)$ holds in $R/\mathfrak{a}$, where $\mathfrak{a} = (a_1,\ldots,a_n)$ and $x \in R^m$ by using induction on the number of reduction rules (R1)–(R4). Inspecting the same reduction rules, one sees that $\text{Red}_{\bar{a}}\,\mathtt{P}(\bar{x})$ has an additional block of existential quantifiers at the rear end of the quantifier string in $\mathtt{P}(\bar{x})$ (when put in prenex-form), from which the last statement now follows immediately. ∎

## 1. The Theory of Artinian Local Rings

**1.1. Proposition.** *The class of all Artinian local rings of length at most $l$ is elementary.*

*Proof.* First, we want to express in first order language that for an element $x$ in a ring $R$, the annihilator $\text{Ann}_R(x)$ is a prime ideal, and whence, by definition an associated prime. The following formula takes care of that,

$$(1.1.1) \qquad \mathtt{Ass}(\boldsymbol{x}) \stackrel{\text{def}}{=} \boldsymbol{x} \neq \boldsymbol{0} \wedge (\forall \boldsymbol{a}, \boldsymbol{b})[\boldsymbol{abx} = \boldsymbol{0} \to (\boldsymbol{ax} = \boldsymbol{0} \vee \boldsymbol{bx} = \boldsymbol{0})].$$

We can also express that the associated prime $\text{Ann}_R(x)$ is maximal by the following formula

$$(1.1.2) \qquad \mathtt{Mass}(\boldsymbol{x}) \stackrel{\text{def}}{=} \mathtt{Ass}(\boldsymbol{x}) \wedge (\forall \boldsymbol{s})[\boldsymbol{sx} \neq \boldsymbol{0} \to (\exists \boldsymbol{r}, \boldsymbol{t})[\boldsymbol{tx} = \boldsymbol{0} \wedge \boldsymbol{t} + \boldsymbol{rs} = \boldsymbol{1}]].$$

Let us now define sentences $\mathtt{Ar}_l$, by induction on $l$, where $l = 1, 2, \ldots$. In case $l = 1$, we define

$$(1.1.3) \qquad \mathtt{Ar}_1 \stackrel{\text{def}}{=} (\forall \boldsymbol{s})[\boldsymbol{s} = \boldsymbol{0} \vee (\exists \boldsymbol{t})[\boldsymbol{st} = \boldsymbol{1}]].$$

In other words, any ring satisfying $\mathtt{Ar}_1$ is a field and therefore of length one (over itself). Define now, using the Reduction Lemma **(0.2)**, by induction

$$(1.1.4) \qquad \mathtt{Ar}_l \stackrel{\text{def}}{=} (\exists \boldsymbol{x})[\mathtt{Mass}(\boldsymbol{x}) \wedge \text{Red}_{\boldsymbol{x}}\,\mathtt{Ar}_{l-1}].$$

Let us prove, again by induction, that a ring satisfying $\mathtt{Ar}_l$ has exactly length $l$ over itself. By the condition $\mathtt{Mass}(\boldsymbol{x})$, which holds for some $x \in R$, we know that $\text{Ann}_R(x)$ is a maximal ideal. Consider the exact sequence

$$(1) \qquad 0 \to xR \longrightarrow R \longrightarrow R/xR \to 0.$$

The first module is isomorphic with $R/\text{Ann}_R(x)$, whence a field and hence of length one over $R$. The last module is of length $l-1$, since $\text{Red}_x\,\mathtt{Ar}_{l-1}$ holds in $R$, so that $\mathtt{Ar}_{l-1}$ holds in $R/xR$, which means, by induction, exactly what we claimed. Therefore it follows from the exact sequence (1) that $R$ is of length $l$. If we add the sentence

$$(1.1.5) \qquad \mathtt{Art}_l = \mathtt{Loc} \wedge \bigvee_{i=1}^{l} \mathtt{Ar}_i$$

to the theory of rings $\mathcal{R}\text{ing}$, we obtain a theory $\mathcal{A}\text{rt}_l$ having as models exactly the Artinian local rings of length at most $l$. ∎

*Remark 1.* From our proof it follows immediately that the class of all Artinian local rings of length exactly $l$ is also elementary. We will denote its theory by $\mathcal{A}\text{rt}'_l$.

*Remark 2.* Note that the class of all local rings of depth zero is given by the single sentence $\mathtt{Loc} \wedge (\exists \boldsymbol{x})[\mathtt{Mass}(\boldsymbol{x})]$, since the latter expresses that a ring is local and its maximal ideal is associated.

**1.2. Observation.** *Let $(R, \mathfrak{m})$ and $(S, \mathfrak{n})$ be Artinian local rings and let $f : R \to S$ be a ring morphism between them. Then $f$ is automatically local.*

*Proof.* Take an element $m \in \mathfrak{m}$. Since $R$ is Artinian, $m$ is nilpotent, and hence also its image $f(m)$ under $f$. Therefore, $f(m)$ cannot be a unit in $S$, whence must lie in $\mathfrak{n}$. ∎

**1.3. Proposition.** *Let $\{R_\alpha, \varphi_{\alpha\beta}\}$ be a direct system of Artinian local rings $R_\alpha$ with maximal ideals $\mathfrak{m}_\alpha$. Let $R = \varinjlim R_\alpha$ be the direct limit. If $l = \limsup \ell(R_\alpha) < \infty$, then $R$ is an Artinian local ring with maximal ideal $\varinjlim \mathfrak{m}_\alpha$, of length $l$.*

*Proof.* Let $\psi_\alpha : R_\alpha \to R$ denote the canonical maps. Set $\mathfrak{m} = \varinjlim \mathfrak{m}_\alpha$. Let us first prove that this is the unique maximal ideal of $R$. So, take any non-unit $r \in R$. Hence, there is some $\Omega$, such that, for all $\alpha \geq \Omega$, we can write $r = \psi_\alpha(r_\alpha)$ for some $r_\alpha \in R_\alpha$. Therefore, $r_\alpha$ is not a unit in $R_\alpha$ and hence must belong to $\mathfrak{m}_\alpha$, so that $r \in \mathfrak{m}$. This proves already that $R$ is local.

Next, we want to show that there exists an $r \in R$, such that $\text{Ann}_R(r) = \mathfrak{m}$. Note that for any Artinian local ring $(S, \mathfrak{n})$ of length at most $e$, we must have that $\mathfrak{n}^e = 0$, i.e., its exponent cannot exceeds its length. Hence, since the lengths of the $R_\alpha$ are bounded for all $\alpha$, say by $e$, we have that $\mathfrak{m}_\alpha^e = 0$. Take now any $e$ elements in $\mathfrak{m}$. We can find a $\Omega$, such that these $e$ elements all are lying in $\psi_\alpha(R_\alpha)$, for all $\alpha \geq \Omega$, and hence in $\psi_\alpha(\mathfrak{m}_\alpha)$. Therefore their product must be zero, since $\mathfrak{m}_\alpha^e = 0$. This proves that $\mathfrak{m}^e = 0$, i.e., $R$ has a finite exponent, say $f$. In particular, $\mathfrak{m}^{f-1} \neq 0$, so that we can take a non-zero element $r$ in it. Since however $\mathfrak{m}^f = 0$, we have that $r\mathfrak{m} = 0$, which gives us the desired element with $\text{Ann}_R(r) = \mathfrak{m}$.

We can now conclude the proof as follows, by induction on $l$. The case $l = 1$ is obvious, since in that case there exists a $\Omega$, such that for all $\alpha \geq \Omega$, the length of $R_\alpha$ is one. For $l > 1$, there is a $\Omega$, such that, for all $\alpha \geq \Omega$, we have that $r = \psi_\alpha(r_\alpha)$ with $r_\alpha \in R_\alpha$ such that $\varphi_{\beta\alpha}(r_\alpha) = r_\beta$, for $\beta \geq \alpha$. Define

$$S_\alpha = \frac{R_\alpha}{r_\alpha R_\alpha}.$$

For $\alpha \geq \Omega$, we have that the length of $S_\alpha$ is at most $l - 1$. It is clear that the $\{S_\alpha\}$ form again a directed system of Artinian local rings with direct limit $R/rR$. Let $h = \limsup \ell(S_\alpha) \leq l - 1$. Since $h < l$, we have by induction, that the length of $R/rR$ is also $h$. By the following exact sequence

$$0 \to rR \longrightarrow R \longrightarrow R/rR \to 0,$$

we therefore obtain that $\ell(R) = h + 1$, since $rR \cong R/\mathfrak{m}$ is of length one. In particular, we proved already that $R$ is Artinian, whence Noetherian. Therefore, the maximal ideal $\mathfrak{m}$ is finitely generated and there exists a $\delta \geq \Omega$, such that

(1) $$\mathfrak{m}_\alpha R = \mathfrak{m},$$

for all $\alpha \geq \delta$. Hence

$$\psi_\alpha(\mathfrak{m}_\alpha r_\alpha) = \mathfrak{m} r = 0$$

in $R$, for all such $\alpha$. Hence there exists an $\beta \geq \alpha$, such that

(2) $$\varphi_{\beta\alpha}(\mathfrak{m}_\alpha r_\alpha) = 0$$

in $R_\beta$. But $\varphi_{\beta\alpha}(\mathfrak{m}_\alpha) = \mathfrak{m}_\beta$ by (1), so that (2) becomes

$$\mathfrak{m}_\beta r_\beta = 0.$$

In other words, we see that the length of $r_\beta R_\beta$ is one. Since we can find for every $\alpha$ a $\beta \geq \alpha$ such that this holds, we conclude that $h = l - 1$. ∎

**1.4. Corollary.** *The class of all Artinian local rings of length at most $l$ (respectively, of length exactly $l$), is $\forall_2$-elementary.*

*Proof.* As seen in **(1.1)**, the models of $\mathcal{A}rt_l$ are exactly the Artinian local rings of length at most $l$. Since $\mathcal{R}ing$ consists out of $\forall_2$-sentences, we proved our claim, if we can prove that $\texttt{Art}_l$ is $\sigma$-persistent with respect to $\mathcal{R}ing$, by using [**Rob**, Theorem 3.4.8].

To prove the $\sigma$-persistency, we only need to show that any increasing chain of Artinian local rings of length at most $l$, has a union with the same property. But this is a consequence of **(1.3)**. The case of $\mathcal{A}rt'_l$ is proved similarly, using $\texttt{Art}'_l = \texttt{Loc} \wedge \texttt{Ar}_l$ instead of $\texttt{Art}_l$. ∎

**1.5. Proposition.** *For each $l$, there exists a $\exists_1$-sentence $\mathtt{Len}_l$, such that, if $(R, \mathfrak{m})$ is an Artinian local ring in which $\mathtt{Len}_l$ holds, then its length is at least $l$. Moreover, if $(R, \mathfrak{m})$ is an Artinian local Gorenstein ring of length at least $l$, then $\mathtt{Len}_l$ holds in it.*

*Remark.* Note that the sentences $\mathtt{Art}_l$ in (1.1.5) are clearly not $\exists_1$.

*Proof.* Let us first give the definition of the $\mathtt{Len}_l$. We proceed by induction on $l$. For $l = 1$, the sentence

$$(1.5.1) \qquad \mathtt{Len}_1 \stackrel{\text{def}}{=} (\exists \boldsymbol{a})[\boldsymbol{a} \neq \boldsymbol{0}]$$

clearly will do. In order to give the recursive description of $\mathtt{Len}_l$, we make the contents of these sentences more explicit. We claim that $\mathtt{Len}_l$ is of the form

$$\mathtt{Len}_l = (\exists \bar{\boldsymbol{a}})(\exists \bar{\boldsymbol{b}})[\chi_l(\bar{\boldsymbol{a}}, \bar{\boldsymbol{b}})],$$

where $\chi_l(\bar{\boldsymbol{a}}, \bar{\boldsymbol{b}})$ is a quantifier free formula with $\mathrm{len}(\bar{\boldsymbol{a}}) = l$ and $\mathrm{len}(\bar{\boldsymbol{b}}) = l - 1$. Moreover, we want that $\chi_l(a, b)$ holds in $R$, with $a_i, b_i \in R$, if and only if, there exists an increasing chain

$$(1) \qquad \mathfrak{a}_0 = 0 \subsetneq \mathfrak{a}_1 \subsetneq \cdots \subsetneq \mathfrak{a}_l$$

in $R$, where each ideal $\mathfrak{a}_i$ is generated by the previous one, $\mathfrak{a}_{i-1}$ and one more element $a_i$. We show how to construct $\mathtt{Len}_{l+1}$, given $\mathtt{Len}_l$ expressing (1). Let

$$(1.5.2) \qquad \mathtt{Len}_{l+1} = (\exists \bar{\boldsymbol{a}})(\exists \bar{\boldsymbol{b}})[\chi_l(\bar{\boldsymbol{a}}, \bar{\boldsymbol{b}}) \wedge \boldsymbol{b_l a_{l+1}} \neq \boldsymbol{0} \wedge \bigwedge_{i=1}^{l} \boldsymbol{b_l a_i} = \boldsymbol{0}].$$

Suppose that $\chi_{l+1}(a, b)$ holds in $R$. To show that this implies the existence of a chain as in (1), of length $l + 1$ this time, we only need to show by induction that the element $a_{l+1}$ is not contained in $\mathfrak{a}_l$, so that the ideal $\mathfrak{a}_{l+1}$ generated by $a_{l+1}$ and $\mathfrak{a}_l$ is strictly bigger than the latter ideal. This, however, is clear, since $b_l$ lies in the annihilator of $\mathfrak{a}_l$, but does not annihilate $a_{l+1}$.

To prove the second statement, assume that $(R, \mathfrak{m})$ is an Artinian local Gorenstein ring. In this case we always have that, if $\mathfrak{a} \subsetneq \mathfrak{b}$ is a strict inclusion of ideals, then $\mathrm{Ann}_R(\mathfrak{b}) \subsetneq \mathrm{Ann}_R(\mathfrak{a})$, since by [**Sch 1**, Theorem 2.2], the operation $\mathrm{Ann}_R(.)$ is an inclusion reversing involution on the lattice of ideals of $R$. If $R$ has length at least $l$, then we can find a chain of ideals $\mathfrak{a}_i = (a_1, \ldots, a_i)$ as in (1). By induction, we show that the $\mathtt{Len}_i$ hold for each $i$. This is clear for $i = 1$, so suppose already established for $i < l$ and we want to show that $\mathtt{Len}_{i+1}$ also holds. By the remark above, the annihilators of $\mathfrak{a}_i$ and $\mathfrak{a}_{i+1}$ must be different, since these ideals are, so that we can find an element $b_i$, annihilating $\mathfrak{a}_i$ but not $a_{i+1}$.   ∎

**1.6. Theorem.** *Let $(R, \mathfrak{m})$ be an existentially closed model of $\mathcal{A}\mathrm{rt}_l$ (or of $\mathcal{A}\mathrm{rt}'_l$). Then the following holds.*
  (i)  *The residue field $\kappa = R/\mathfrak{m}$ is algebraically closed.*
  (ii) *$R$ is Gorenstein.*
  (iii) *The length of $R$ is exactly $l$.*

*Proof.* (i). Let $p$ be a monic irreducible polynomial in one variable $T$ over $\kappa$. Take a monic lifting $P$ of $p$ to $R$ and consider

$$\tilde{R} = \frac{R[T]}{(P)}.$$

Let us prove that $\tilde{R}$ is local with maximal ideal $\tilde{\mathfrak{m}} = \mathfrak{m}\tilde{R}$. First, $\tilde{\mathfrak{m}}$ is a maximal ideal since

$$\tilde{R}/\mathfrak{m}\tilde{R} \cong \frac{\kappa[T]}{(p)},$$

which is a field since $p$ is irreducible over $\kappa$. Moreover, $\tilde{R}$ is a finite free $R$-algebra, and hence in particular Artinian and whence semi-local. Therefore, each maximal ideal $\mathfrak{n}$ of $\tilde{R}$, must contract to $\mathfrak{m}$ in $R$, and hence contain $\tilde{\mathfrak{m}}$, whence be equal to it. This proves that $\tilde{R}$ is local.

Let us prove, by induction on $l$, that the length of $\tilde{R}$ is at most $l$, so that $\tilde{R}$ is also a model of $\mathcal{A}\mathfrak{rt}_l$. (To this end we will not yet use the existentially closedness of $R$). The case $l = 1$ is clear, since then also $\tilde{R}$ a field. Assume therefore that $l > 1$ and choose $x \in R$, such that $\mathrm{Ann}_R(x) = \mathfrak{m}$. Then clearly also $\mathrm{Ann}_{\tilde{R}}(x) = \tilde{\mathfrak{m}}$, so that $x\tilde{R}$ is of length one. But

$$\tilde{R}/x\tilde{R} = \frac{(R/xR)[T]}{(\bar{P})},$$

where $\bar{P}$ denotes the image of $P$ in $(R/xR)[T]$, which also reduces to $p$ modulo the maximal ideal of $R/xR$. Clearly $R/xR$ is a model of $\mathcal{A}\mathfrak{rt}_{l-1}$. We can apply the induction hypothesis to $\tilde{R}/x\tilde{R}$, to conclude that its length is at most $l - 1$. Our claim then follows from the exact sequence

$$0 \to x\tilde{R} \longrightarrow \tilde{R} \longrightarrow \tilde{R}/x\tilde{R} \to 0$$

and the observation that $\ell(x\tilde{R}) = 1$.

Hence, $\tilde{R}$ is a model of $\mathcal{A}\mathfrak{rt}_l$, admitting $T$ as a root of $P$. By existentially closedness of $R$, we must have that there exists already a root $r$ of $P$ in $R$, whence its residue modulo $\mathfrak{m}$ is a root of $p$, proving that $\kappa$ is algebraically closed.

(*ii*). Let $0 \neq \eta \in \mathrm{Soc}\, R$. We want to show that $\eta$ generates the socle, showing that $R$ is Gorenstein. Suppose not, so that there exists $r \in \mathrm{Soc}\, R$ with $\eta \notin (r)$. Let $T$ be one variable and let $\Sigma$ be the multiplicative set of $R[T]$ consisting of all monic polynomials. Let $R(T)$ denote the localisation $\Sigma^{-1}R[T]$. The ideal $\mathfrak{m}R(T)$ is maximal, since the quotient equals $\kappa(T)$, which is a field. It is moreover nilpotent and hence it must be the unique maximal ideal of $R(T)$, showing that the latter is an Artinian local ring. By induction on the length, one shows that $R(T)$ has the same length as $R$.[3] Let $\tilde{R} = R(T)/(\eta + rT)$, which therefore is also a model of $\mathcal{A}\mathfrak{rt}_l$. We claim that $R \subset \tilde{R}$. Assuming the claim and using that $R$ is existentially closed, we see that the the equation $\eta + rT = 0$ viewed as an equation in the variable $T$ has a root in $\tilde{R}$ and whence in $R$, contradicting that $\eta \notin (R)$.

So let us show the claim. Let $c \in R$, such that $c = 0$ in $\tilde{R}$. This means that there exist $a_i, b_i \in R$, such that

$$c(b_0 + b_1 T + \cdots + b_s T^s + T^{s+1}) = (\eta + rT)(a_0 + a_1 T + \cdots + a_s T^s)$$

in $R[T]$. This leads to a system of equalities in the coefficients

$$c = ra_s$$
$$cb_s = ra_{s-1} + a_s\eta$$
$$cb_{s-1} = ra_{s-2} + a_{s-1}\eta$$
$$\vdots$$
$$cb_1 = ra_0 + a_1\eta$$
$$cb_0 = a_0\eta.$$

The first of these equalities tells us that $(c) \subset (r)$, hence the remaining imply that $a_i\eta \in (r)$. By our assumption, none of the $a_i$ can therefore be a unit, so must lie in $\mathfrak{m}$. Hence from the first equality it then follows that $c = 0$, since $r \in \mathrm{Soc}\, R$.

(*iii*). Suppose not, so that, say, $\ell(R) = h < l$. Let

$$\tilde{R} = \frac{R[T]}{(T^2, \mathfrak{m}T)}.$$

As can be easily checked, $\tilde{R}$ is again an Artinian local ring with maximal ideal $\tilde{\mathfrak{m}} = (\mathfrak{m}, T)$. One also checks that $\mathrm{Ann}_{\tilde{R}}(T) = \tilde{\mathfrak{m}}$, so that, by considering the exact sequence

$$0 \to (T) \longrightarrow \tilde{R} \longrightarrow R \to 0$$

---

[3] By induction $R(T)/\eta R(T) \cong (R/\eta R)(T)$ has length $\ell(R) - 1$. But the maximal ideal of $R(T)$ annihilates $\eta R(T)$, so the latter ideal has length one and we are done.

we deduce that $\tilde{R}$ is of length $h + 1 \leq l$. By **(1.4)** the theory $\mathcal{A}\mathfrak{rt}'_{h+1}$ of all Artinian local rings of length exactly $h + 1$ is $\forall_2$-elementary, so that by [**Hod 2**, Corollary 8.2.3], we can embed $\tilde{R}$ in an existentially closed model $S$ of that theory. By $(i)$ above, we know that $S$ is Gorenstein and clearly of length $h + 1$ (since all models of this theory have fixed length). By **(1.5)**, the $\exists_1$-sentence $\texttt{Len}_{h+1}$ holds in $S$. Since, however, $S$ is also a model of $\mathcal{A}\mathfrak{rt}_l$, by assumption on $h$, and since $S$ is an extension of $R$, which is assumed to be existentially closed for $\mathcal{A}\mathfrak{rt}_l$, we conclude that also $\texttt{Len}_{h+1}$ holds in $R$. But by **(1.5)** again, this means that $R$ has length at least $h + 1$, contradiction. ∎

**1.7. Corollary.** *Let $(R, \mathfrak{m})$ be an Artinian local ring of length $l$. Then there exists an Artinian local Gorenstein ring $(S, \mathfrak{n})$ of the same length $l$ with an algebraically closed residue field, such that $R \subset S$.*

*Proof.* Since the theory $\mathcal{A}\mathfrak{rt}'_l$ is $\forall_2$-elementary, every model can be embedded in an existentially closed model of the theory (see for instance [**Hod 2**, Theorem 8.2.1]). Together with **(1.6)**, this provides the proof of our claim. ∎

**1.8. Example.** Let $\kappa$ be a field and let $R = \kappa[S,T]/\mathfrak{a}$, where $\mathfrak{a} = (S^2, ST^2, T^3)$. One easily verifies that $R$ is an Artinian local ring of length 5, embedding dimension 2, exponent 3 and type 2. In particular, $R$ is not Gorenstein. The following Artinian local ring $\tilde{R}$ provides an extension of $R$ of the same length which is moreover Gorenstein. Let $\lambda$ be the fraction field of $\kappa[A,B,C]/(A^2+B^2+C^2)$ and let $\tilde{R} = \lambda[X,Y,Z]/\mathfrak{b}$, where $\mathfrak{b} = (X^2-Y^2, X^2-Z^2, XY, XZ, YZ)$. The reader should check that this is indeed an Artinian local Gorenstein ring of length 5 (embedding dimension 3 and exponent 3). The embedding $\varphi: R \to \tilde{R}$ is now given by sending (the image of) $S$ to $AX + BY + CZ$ and sending $T$ to $BX$. Since the exponent of $\tilde{R}$ is 3, it is clear that both $ST^2$ and $T^3$ are sent to zero under $\varphi$. The image of $S^2$ under $\varphi$ equals $(AX + BY + CZ)^2$. But the latter is equal to $(A^2 + B^2 + C^2)X^2$ in $\tilde{R}$. Moreover $A^2 + B^2 + C^2 = 0$ in $\lambda$, hence also in $\tilde{R}$. This proves that the map $\varphi$ is well defined. To check that it is injective, let us calculate the image of an arbitrary element $f$ of $R$. Each element of $R$ can be (uniquely) written as

$$f = u + vS + wT + xST + yT^2,$$

where $u, v, w, x, y$ are elements of $\kappa$. Hence

$$\varphi(f) = u + v(AX + BY + CZ) + wBX + xABX^2 + yB^2X^2.$$

For this element to be zero, we must have that $u = 0$, $vA + wB = 0$, $vB = 0$, $vC = 0$ and $xAB + yB^2 = 0$. Hence $u = v = w = 0$ and $xA + yB = 0$. However, since $A$ and $B$ are transcendent over $\kappa$, the latter equality also implies that $x = y = 0$, proving that $f = 0$ in $R$.

Note that we didn't need the residue field of $\tilde{R}$ to be algebraically closed, but we did need it to be transcendent over $\kappa$! Note also that $\tilde{R}$ is not a complete intersection ring.

## 2. The Equicharacteristic Case

**2.1. Notation.** In this and the next section, $(R, \mathfrak{m})$ will be an Artinian local ring with residue field $\kappa$. We fix once and for all a set of minimal generators $x = (x_1, \ldots, x_\mu)$ for the maximal ideal $\mathfrak{m}$. Also, on the set of $\mu$-tuples $\mathbb{N}^\mu$, we choose once and for all a *monomial ordering* $\leq$. This means that $\leq$ totally orders $\mathbb{N}^\mu$ with 0 as the minimal element such that for any $\alpha, \beta, \Omega \in \mathbb{N}^\mu$, we have that $\alpha \leq \beta$ implies that $\alpha + \Omega \leq \beta + \Omega$. For future purposes, we might also want to require that $\leq$ is *graded*, with which we mean that whenever $|\alpha| < |\beta|$, then $\alpha < \beta$, where we write $|(a_1, \ldots, a_\mu)| = a_1 + \cdots + a_\mu$. In other words, $\leq$ respects the grading on $\mathbb{N}^\mu$. An example of a graded monomial ordering is the so called *graded lexicographical* ordering, which is defined by calling $\alpha$ smaller than $\beta$, if $|\alpha| < |\beta|$, or, if $|\alpha| = |\beta|$, but $\alpha$ is lexicographically smaller than $\beta$. In this paper, however, we can take the lexicographical order on $\mathbb{N}^\mu$, or any other monomial ordering, for that matter. Therefore, for sake of exposition, we will just assume that $\leq$ is the usual lexicographical ordering.

We define the *support* $\text{Supp}(R)$ as the set of all $\mu$-tuples $\alpha$ in $\mathbb{N}^\mu$, such that $x^\alpha \neq 0$ in $R$. Note that this is a finite set. We define, for any $\alpha \in \text{Supp}(R)$, the ideal $\mathfrak{a}_R(\alpha)$ (for short, $\mathfrak{a}(\alpha)$) as the ideal generated by all $x^\beta$, with $\beta > \alpha$. We will use the notation $\Delta(R, x, \leq)$ (or just $\Delta_R$, when the choice of a minimal set of generators $x$ and the choice of a monomial ordering is understood) for the collection of all $\alpha$, for which $x^\alpha \notin \mathfrak{a}(\alpha)$. We denote the collection of all $x^\alpha$ with $\alpha \in \Delta_R$ by $\mathcal{E}_R$.

In this section we furthermore assume that $R$ has the same characteristic as $\kappa$, the so called *equicharacteristic* case. Under this additional assumption, $R$ is a finite dimensional local $\kappa$-algebra, generated by the $x_i$, and conversely, any such ring is of the above type, by COHEN's Structure Theorem. By definition of monomial ordering, any $\beta + \Omega \geq \beta$, so that $\mathfrak{a}(\alpha)$ equals the vector space over $\kappa$ spanned by all the $x^\beta$, with $\beta \in \mathrm{Supp}(R)$ and $\beta > \alpha$. In other words, if $\alpha \notin \Delta_R$, then $x^\alpha$ equals some $\kappa$-linear combination of the $x^\beta$, with $\beta > \alpha$.

**2.2. Lemma.** *Let $R$ be an equicharacteristic Artinian local ring with residue field $\kappa$ and let $\Delta_R = \Delta(R, x, \leq)$ be as defined in* **(2.1)**. *Then the set $\mathcal{E}_R = \{ x^\alpha \mid \alpha \in \Delta_R \}$ is a base of $R$ over $\kappa$.*

*Proof.* Let us first prove that the $\kappa$-vector space $V$ generated by $\mathcal{E}_R$ is $R$ itself. If not, then there exists a maximal $\alpha \in \mathrm{Supp}(R)$, such that $x^\alpha \notin V$. In particular $\alpha \notin \Delta_R$, hence, there exist $r_\beta \in \kappa$ (see remark above), such that

(1) $$x^\alpha = \sum_{\beta > \alpha} r_\beta x^\beta.$$

By the maximality of $\alpha$, all $x^\beta \in V$, so that by (1), also $x^\alpha \in V$, contradiction.

Let us now show that the $x^\Omega$, for $\Omega \in \Delta_R$ are $\kappa$-linearly independent. Suppose not, then there exists a $\Gamma \subset \Delta_R$ and non-zero $u_\Omega \in \kappa$, for $\Omega \in \Gamma$, such that

(2) $$\sum_{\Omega \in \Gamma} u_\Omega x^\Omega = 0.$$

Let $\alpha$ be the smallest element of $\Gamma$, so that by (2), we have that $x^\alpha \in \mathfrak{a}(\alpha)$, contradicting that $\alpha \in \Delta_R$. ∎

**2.3. Proposition.** *Let $(R, \mathfrak{m})$ and $(S, \mathfrak{n})$ be equicharacteristic Artinian local rings. Suppose that $R \subset S$ and that the length of $S$ does not exceed the length of $R$. Let $\kappa$ (respectively, $\lambda$) be the residue field of $R$ (respectively, $S$). If $R$ is Gorenstein, then $S = R \otimes_\kappa \lambda$ (for some embedding of $\kappa$ in $\lambda$).*

*Proof.* By COHEN's Structure Theorem for complete local rings, we can write $R = \kappa[X]/T$, for some $(X)$-primary ideal $I$. Let $x = (x_1, \ldots, x_n)$ denote the image of $X = (X_1, \ldots, X_n)$ in $R$. Let $l$ be the length of $R$. We adopt the notations of **(2.1)** for this ring $R$. In particular, we will write $\Delta_R$ for $\Delta(R, x, \leq)$. Since $\kappa \subset S$, we can find by [**Mats**, Theorem 28.3] a field $\kappa \subset \lambda' \subset S$ which is isomorphic with the residue field $\lambda$ of $S$. Hence without loss of generality, we take $\lambda = \lambda'$, so that the embedding $R \subset S$ is in fact a $\kappa$-algebra morphism.

We claim that the basis $\mathcal{E}_R = \{ x^\alpha \mid \alpha \in \Delta_R \}$ remains independent over $\lambda$ in $S$ (where we still write $x$ for its image in $S$). Assuming the claim, then since $S$ has dimension at most $l$ (over $\lambda$), $\mathcal{E}_R$ must be a base of $S$ over $\lambda$ as well. In particular, the canonical map

$$R \otimes_\kappa \lambda \to S : r \otimes u \mapsto ur$$

is surjective. Since both sides however have the same dimension over $\lambda$, the map must be an isomorphism.

So we only need to show our claim. Hence suppose the contrary, so that there exist $u_\Omega \in \lambda$, not all zero, for which

(1) $$\sum_{\Omega \in \Delta_R} u_\Omega x^\Omega = 0$$

in $S$. Let $\Gamma = \{ \Omega \in \Delta_R \mid u_\Omega \neq 0 \}$ Let $\alpha \in \Gamma$ be minimal in $\Gamma$. Let $\mathfrak{a}$ be the ideal of $R$ generated by all the $x^\Omega$ with $\Omega \in \Gamma$ different from $\alpha$. For $\Omega \in \Gamma$ different from $\alpha$, we have by minimality of $\alpha$ that $\alpha < \Omega$. Hence $x^\Omega \in \mathfrak{a}(\alpha)$. Therefore, we have that $\mathfrak{a} \subset \mathfrak{a}(\alpha)$. Since $\alpha \in \Delta_R$, we have that $x^\alpha \notin \mathfrak{a}(\alpha)$. Hence, in particular, we have that

(2) $$x^\alpha \notin \mathfrak{a}.$$

Since $R$ is Gorenstein, we have by [**Sch 1**, Theorem 2.2] that $\operatorname{Ann}_R(\cdot)$ is an inclusion reversing involution on the lattice of ideals of $R$. Hence, from (2), it follows that $\operatorname{Ann}_R(\mathfrak{a}) \not\subset \operatorname{Ann}_R(x^\alpha)$. Let $s \in \operatorname{Ann}_R(\mathfrak{a})$, such that $sx^\alpha \neq 0$. Since $R$ embeds in $S$, the latter also holds in $S$. But this leads to a contradiction when we multiply (1) with $s$. This therefore shows the validity of our claim and hence of the proposition. ∎

*Remark.* Note that this implies that $S$ is faithfully flat over $R$ (since $\lambda$ is faithfully flat over $\kappa$) and that $S$ is also Gorenstein. Note also that in case $S$ has length bigger than $l$, we still have an inclusion $R \otimes_\kappa \lambda \hookrightarrow S$.

**2.4. Proposition.** *Let $(R, \mathfrak{m})$ be an equicharacteristic Artinian local ring with an algebraically closed residue field $\kappa$. Let $\lambda$ be an extension field of $\kappa$ and set $S = R \otimes_\kappa \lambda$. Let $P_i(T) \in R[T]$ be polynomials, for $i = 1, \ldots, N$, where $T = (T_1, \ldots, T_n)$ are variables. Suppose there exists a $s \in S^n$, such that $P_i(s) = 0$, for $i = 1, \ldots, M$, and $P_i(s) \neq 0$, for $i = M+1, \ldots, N$. Then we can find a $r \in R^n$, satisfying the same equalities and inequalities, i.e., $P_i(r) = 0$, for $i = 1, \ldots, M$ and $P_i(r) \neq 0$, for $i = M+1, \ldots, N$.*

*Proof.* Since $R$ is equicharacteristic, it is of the form $\kappa[X]/I$, for some variables $X$ and some $(X)$-primary ideal $I$. Let us write $P_i(T, X)$ for the liftings of the $P_i$ to $\kappa[T, X]$. Let $\Delta_R$ be as before, such that the (images of the) $X^\delta$, for $\delta \in \Delta_R$, form a basis of $R$ over $\kappa$. Clearly, the same monomials give a basis of $S$ over $\lambda$.

For $j = 1, \ldots, n$, let
$$f_j(A, X) = \sum_{\delta \in \Delta_R} A_{j\delta} X^\delta,$$
where $A$ is the collection of all variables $A_{j\delta}$ and let $f = (f_1, \ldots, f_n)$. Hence, for each $i = 1, \ldots, N$, we can find polynomials $p_{i\delta}(A) \in \kappa[A]$, such that

$$P_i(f, X) = \sum_{\delta \in \Delta_R} p_{i\delta}(A) X^\delta \tag{1}$$

over $R$. Finally, we can find $b_{j\delta} \in \lambda$, for $j = 1, \ldots, n$ and $\delta \in \Delta_R$, such that

$$f(b, X) = s \tag{2}$$

in $S$, where $b = (b_{j\delta})$. Using (1) and (2), we obtain, for all $\delta \in \Delta_R$ and all $i = 1, \ldots, M$, that $p_{i\delta}(b) = 0$, and, moreover, for each $i = M+1, \ldots, N$, there exists an $\delta_i \in \Delta_R$, such that $p_{i\delta_i}(b) \neq 0$. Since $\kappa$ is algebraically closed (and hence existentially closed for the theory of fields), we can find $a_{i\delta} \in \kappa$, such that $a = (a_{i\delta})$ satisfies the same equalities and inequalities above, upon replacing $a$ by $a$. I.e., for all $\delta \in \Delta_R$ and all $i = 1, \ldots, M$, we have that $p_{i\delta}(a) = 0$, and, for all $i = M+1, \ldots, N$, that $p_{i\delta_i}(a) \neq 0$. Therefore, if we put $f(a, X) = \bar{r}$, then this $n$-tuple over $R$ fulfills our requirements. ∎

**2.5. Theorem.** *The equicharacteristic existentially closed models of $\mathcal{A}\mathrm{rt}_l$ are exactly the equicharacteristic Artinian local Gorenstein rings of length $l$ with algebraically closed residue field.*

*Proof.* From **(1.6)** we know already that existentially closed models of $\mathcal{A}\mathrm{rt}_l$ are Gorenstein, of length $l$ and their residue field is algebraically closed.

Conversely, let $(R, \mathfrak{m})$ be an equicharacteristic Artinian local Gorenstein ring of length $l$ with algebraically closed residue field $\kappa$. Let $(S, \mathfrak{n})$ be a model of $\mathcal{A}\mathrm{rt}_l$ with residue field $\lambda$, such that $R \subset S$ (so that in particular also $S$ is equicharacteristic). Since the length of $S$ is at most $l$, we know from **(2.3)** that $S \cong R \otimes_\kappa \lambda$. We have to show that any system of equalities and inequalities over $R$ which has a solution over $S$, has already a solution over $R$. This is the contents of **(2.4)**. ∎

3. THE MIXED CHARACTERISTIC CASE

**3.1. Definition.** In this section we take up the same conventions and notations as in **(2.1)**. In particular $(R, \mathfrak{m}, \kappa)$ is an Artinian local ring. However, this time we assume that $\kappa$ has positive characteristic $p$, but $p \neq 0$ in $R$. This means that $\kappa$ is not embeddable in $R$. Nevertheless, some power of $p$ must vanish, since any element in $\mathfrak{m}$ is nilpotent. As we will only be interested in the case that the residue field is algebraically closed, we make this assumption throughout as well. Let $e$ be the exponent of $R$.

We claim that there exists some power $q$ of $p$, such that $q = 0$ in $R$ and, for any two elements $a, b \in R$, we have that $a \equiv b \mod \mathfrak{m}$ implies that $a^q = b^q$. Indeed, write $a = b + m$, for some $m \in \mathfrak{m}$. By NEWTON's binomial expansion

$$(1) \qquad a^q = (b+m)^q = b^q + \binom{q}{1} b^{q-1} m + \binom{q}{2} b^{q-2} m^2 + \cdots + \binom{q}{e-1} b^{q-e+1} m^{e-1}.$$

Since only the $e$ first terms in (1) survive, we can choose $q$ big with respect to $e$, so that any binomial $\binom{q}{i}$ appearing in (1) is divisible by $p^e$ and hence vanishes in $R$. Note that by the same reasoning, if $a^q = b^q$, we also must have that $a \equiv b \mod \mathfrak{m}$. We fix this $q$ for the rest of the section as well. The reader should note that any power of $p$ big enough yields the same results, so it will not be essential what power $q$ we take.

We are now able to define the *Witt map* $\omega \colon \kappa \to R$ as follows. Let $u \in \kappa$. Since we assumed that $\kappa$ is algebraically closed (in fact, being perfect suffices), we can find a $v \in \kappa$, such that $u = v^q$. Let $b \in R$ such that its residue $\bar{b}$ equals $v$. Then we set $\omega(u) = b^q$. By our choice of $q$, this is independent of the particular choice of the lifting $b$. Moreover, one easily verifies that the residue of $\omega(u)$ is $\bar{b}^q = u$ again and that $\omega$ is a multiplicative map. However, it is not necessarily additive.[4]

The image $\operatorname{Im}\omega$ of the Witt map equals the collection $W = \{a^q \mid a \in R\}$ of *Witt vectors* of $R$. Note that if we were to take even a larger power of $p$ than $q$, the set $W$ would not change in view of the algebraic closedness of $\kappa$ and whence the Witt vectors are independent of the particular choice of $q$, i.e., the Witt vectors are an invariant of $R$. Any non-zero Witt vector is necessarily a unit (if $a^q \in \mathfrak{m}$, then also $a \in \mathfrak{m}$ and hence $a^q = 0$). We have that $R = W + \mathfrak{m}$.

We will call a subset $\mathcal{B} = \{b_1, \ldots, b_k\}$ of $R$ a *Witt base* of $R$, if for any element $r \in R$, there exist uniquely determined $u_i \in \kappa$, such that

$$(1) \qquad r = \sum_{i=1}^{d} \omega(u_i) b_i.$$

The following lemma is elementary but crucial.

**Lemma 3.2.** *Let $w_1, w_2$ be two Witt vectors. Then $w_1 \equiv w_2 \mod \mathfrak{m}$, if and only if, $w_1 = w_2$.*

*Proof.* Suppose $w_1 \equiv w_2 \mod \mathfrak{m}$. Hence by our choice of $q$, we must have that $w_1^q = w_2^q$. Let $w_i = a_i^q$, so that $a_1^{q^2} = a_2^{q^2}$. As observed before, replacing $q$ by $q^2$, this implies that $a_1 \equiv a_2 \mod \mathfrak{m}$ and whence $w_1 = a_1^q = a_2^q = w_2$. ∎

*Remark.* This shows that the decomposition of an element $r \in R$ as $w + m$, with $w \in W$ and $m \in \mathfrak{m}$ is unique.

**3.3. Lemma.** *Let $(R, \mathfrak{m})$ be an Artinian local ring of mixed characteristic, with algebraically closed residue field $\kappa$. Let notation be as in **(2.1)** and **(3.1)**. Then $\mathcal{E}_R$ is a Witt base for $R$.*

*Proof.* Let us first show that for each $r \in R$, there exist $a_\alpha$, for $\alpha \in \Delta_R$, such that

$$(1) \qquad r = \sum_{\alpha \in \Delta_R} a_\alpha^q x^\alpha.$$

---

[4]In fact, it will only be additive provided $p = 0$ in $R$ (in order for all binomials $\binom{q}{n}$ to vanish) and in this case we retrieve $\kappa$ as COHEN's coefficient field for the equicharacteristic Artinian (whence complete) local ring $R$.

Let $r \in R$ be non-zero and write $v(r) = k$ to mean that $r \in \mathfrak{m}^k \setminus \mathfrak{m}^{k+1}$. Formally, put $v(0) = e + 1$. Let us prove by downwards induction on $k$ that any $r \in R$ with $v(r) = k$ can be written in the form (1) with all $a_\alpha = 0$ for $|\alpha| < k$. Note that $r = 0$ trivially satisfies these assumptions. Hence let $v(r) = k$ and the claim proven for all $s \in R$ with $v(s) > k$.

We can write

$$r = \sum_{|\delta|=k} r_\delta x^\delta \tag{2}$$

with $r_\delta \in R$. Each $\delta \notin \Delta_R$ with $|\delta| = k$, can be written in terms of the $x^\alpha$ for $\alpha \in \Delta_R$ with $|\alpha| = k$ and some elements in $\mathfrak{m}^{k+1}$. Let $\Delta_k$ be the set of all indices $\alpha \in \Delta_R$ for which $|\alpha| = k$. After this simplification (2) can be rewritten as

$$r = \sum_{\alpha \in \Delta_k} s_\alpha x^\alpha + r_1,$$

where $r_1 \in \mathfrak{m}^{k+1}$ and $s_\alpha \in R$.

Since $R = W + \mathfrak{m}$, we can find $a_\alpha \in R$ and $m_\alpha \in \mathfrak{m}$, such that $s_\alpha = a_\alpha^q + m_\alpha$, for all $\alpha \in \Delta_k$. This yields an $r_2 \in \mathfrak{m}^{k+1}$, such that

$$r = \sum_{\alpha \in \Delta_k} a_\alpha^q x^\alpha + r_2. \tag{3}$$

By the induction hypothesis, $r_2$ can be written as a linear combination of the $x^\delta$, for $\delta \in \Delta_R$ with $|\delta| > k$, with coefficients in $W$. Together with (3), this shows our claim for $v(r) = k$.

Let us now show that the representation (1) is also unique. For if it were not, then we could find $b_\alpha \in R$, such that

$$r = \sum_{\alpha \in \Delta_R} b_\alpha^q x^\alpha.$$

Subtracting this from (1) yields

$$0 = \sum_{\alpha \in \Delta_R} (a_\alpha^q - b_\alpha^q) x^\alpha. \tag{4}$$

Let $\Omega$ be the lexicographically smallest index for which $0 \neq a_\Omega^q - b_\Omega^q$. If the latter difference were in $\mathfrak{m}$, then by **(3.2)** it would actually be zero. So $a_\Omega^q - b_\Omega^q$ must be a unit, exhibiting that $x^\Omega \in \mathfrak{a}(\Omega)$, contradiction. ∎

*Remark.* Let $\Omega \in \Delta_R$ and suppose $r \in \mathfrak{a}(\Omega)$. Write

$$r = \sum_{\alpha \in \Delta_R} a_\alpha^q x^\alpha.$$

By the same argument used to prove (1), we can show that then all $a_\alpha^q = 0$, for $\alpha \leq \Omega$. Hence, if $\delta \in \Delta_R$ is the smallest index in $\Delta_R$ strictly larger than $\Omega$, then $r \equiv a_\delta^q x^\delta \mod \mathfrak{a}(\delta)$. It follows that $\mathfrak{a}(\Omega)/\mathfrak{a}(\delta) \cong \kappa$. Hence the $\mathfrak{a}(\alpha)$ form a maximal ascending chain of ideals of $R$, when $\alpha$ runs over $\Delta_R$ from the (lexicographically) largest index to the smallest. In particular, the cardinality of $\Delta_R$ equals the length $l$ of the ring $R$.

**3.4. Proposition.** *Let $(R, \mathfrak{m})$ be an Artinian local ring of mixed characteristic, with algebraically closed residue field $\kappa$. Let notation be as in **(2.1)** and **(3.1)**. Let $S$ be an overring of $R$ of the same length $l$. If $R$ is Gorenstein, then $S$ is an Artinian local Gorenstein ring with maximal ideal $\mathfrak{m}S$ and $\mathcal{E}_R$ is a Witt base of $S$. Moreover, $S$ is faithfully flat over $R$.*

*Proof.* Let $0 = \alpha_1 < \alpha_2 < \cdots < \alpha_l$ be the $l$ elements of $\Delta_R$. From the remark after **(3.3)**, we know that

$$0 = \mathfrak{a}(\alpha_l) \subsetneq \mathfrak{a}(\alpha_{l-1}) \subsetneq \cdots \subsetneq \mathfrak{a}(\alpha_1) = \mathfrak{m}, \tag{1}$$

is a maximal chain of ideals in $R$. Let us write $\eta = x^{\alpha_l}$, so that if $R$ is Gorenstein, then $\mathfrak{a}(\alpha_{l-1}) = \eta R$ is the unique non-zero minimal ideal, equal to the socle $\mathrm{Soc}(R)$ of $R$. Moreover, each ideal is an annihilator ideal, so the extension of the chain (1) to $S$ remains strict. In particular, since $S$ has length $l$, we must have that $\mathfrak{m}S$ is a maximal ideal of $S$. Since it is clearly nilpotent, it is the unique maximal ideal of $S$ and hence $S$ is local with maximal ideal $\mathfrak{m}S$. Let $\lambda$ denote its residue field, so that $S \otimes_R \kappa = \lambda$.

Let us define similarly $\Delta_S$ as the collection of all indices $\alpha$ for which $x^\alpha$ does not lie in the ideal $\mathfrak{a}_S(\alpha) = \mathfrak{a}(\alpha)S$ of $S$ generated by all $x^\beta$ with $\beta > \alpha$. Since $x^\alpha \notin \mathfrak{a}(\alpha)$, for $\alpha \in \Delta_R$, there exists $r_\alpha \in R$, such that $r_\alpha x^\alpha \neq 0$ but $r_\alpha \mathfrak{a}(\alpha) = 0$. In particular, we cannot have that $x^\alpha \in \mathfrak{a}(\alpha)S$, so that $\alpha \in \Delta_S$. But clearly $\Delta_S \subset \Delta_R$ and hence we have equality, showing that $\mathcal{E}_R$ is a Witt base of $S$.

Since $r_\alpha x^\alpha \neq 0$, for $\alpha \in \Delta_R$, we can find $r'_\alpha \in R$, such that $r'_\alpha r_\alpha x^\alpha = \eta$. Put $\tau_\alpha = r'_\alpha r_\alpha$, so that $\tau_\alpha x^\alpha = \eta$ and $\tau_\alpha \mathfrak{a}(\alpha) = 0$. Note that $\tau_\alpha \in \mathfrak{m}$ unless $\alpha = \alpha_l$. Let $s \in \mathrm{Soc}(S)$ and write $s$ as

$$(2) \qquad s = \sum_{\alpha \in \Delta_R} s_\alpha^q x^\alpha,$$

for some $s_\alpha \in S$. Let $\Omega \in \Delta_R$ be the smallest element for which $s_\Omega^q$ is non-zero (whence a unit). If $\Omega < \alpha_l$, then multiplying (2) with $\tau_\Omega \in \mathfrak{m}$ yields $0 = \tau_\Omega s_\Omega^q x^\Omega = s_\Omega^q \eta$, contradiction. Hence $\mathrm{Soc}(S)$ is contained in $\eta W \subset \eta S$ and whence equal to the latter ideal, showing that $S$ is Gorenstein.

Let us now proof the last statement. Let $R_i = R/\mathfrak{a}(\alpha_i)$ and $S = S/\mathfrak{a}(\alpha_i)S$. It follows that $\ell(R_i) = \ell(S_i) = i$ and that $S_i = S \otimes_R R_i$. In particular, $R_1 \cong \kappa$ and $S_1 \cong \lambda$. We claim that $\mathrm{Tor}_1^R(S, R_i) = 0$, for all $i$. Assuming the claim, we have, by taking $i = 1$, that $\mathrm{Tor}_1^R(S, \kappa) = 0$, so that $S$ is (faithfully) flat over $R$ by the Local Flatness Criterion [**Mats**, Theorem 22.3]. Hence remains to show our claim, which we will do by downwards induction on $i$. If $i = l$, then $R_l = R$ and the claim trivially holds. Hence assume $i < l$. By the maximality of the chain (1), we have that $\mathfrak{a}(\alpha_{i+1})R_i$ has length one, whence is isomorphic with $\kappa$. Therefore we have an exact sequence

$$0 \to \kappa \longrightarrow R_i \longrightarrow R_{i+1} \to 0.$$

Tensoring with $S$ and our induction hypothesis yields an exact sequence

$$0 = \mathrm{Tor}_1^R(S, R_i) \to \mathrm{Tor}_1^R(S, R_{i+1}) \to S \otimes_R \kappa = \lambda \xrightarrow{\imath} S_i \to S_{i+1} \to 0.$$

Since $S_i$ and $S_{i+1}$ have different lengths, the map $\imath$ cannot be zero and whence must be injective, showing that $\mathrm{Tor}_1^R(S, R_{i+1}) = 0$, as required. ∎

**3.5. Lemma.** *Let $(R, \mathfrak{m})$ be an Artinian local ring of mixed characteristic, with algebraically closed residue field $\kappa$. Let notation be as in (2.1) and (3.1). Let $Q \in R[T]$, where $T$ is a finite set of variables $(T_1, \ldots, T_n)$. Then there exist $h_i \in R[T]$, for $i = 1, \ldots, \mu$ with the following property. For any $v = (v_1, \ldots, v_n) \in \kappa^n$, we have an equality*

$$Q(\omega(v^q)) = \omega(\bar{Q}(v^q)) + \sum_{i=1}^{\mu} x_i h_i(\omega(v)),$$

*where $v^q$ is shortcut for $(v_1^q, \ldots, v_n^q)$ and $\omega(v)$ for $\omega(v_1, \ldots, v_n)$.*

*Proof.* Let $c_i \in R$ be liftings of the $v_i$ and set $c = (c_1, \ldots, c_n)$. In particular, we have (in vector notation) that $\omega(v^q) = c^q$. Writing out the coefficients of $Q$ in the Witt base $\mathcal{E}_R$, we can write

$$Q = \sum_{\alpha \in \Delta_R} \sum_\nu r_{\alpha,\nu}^q T^\nu x^\alpha,$$

with $r_{\alpha,\nu} \in R$. Let $P$ denote the polynomial

$$(1) \qquad P = \sum_{\alpha \in \Delta_R} \sum_\nu r_{\alpha,\nu} T^\nu x^\alpha.$$

One checks that

$$(2) \qquad P(T)^q \equiv Q(T^q) \mod \mathfrak{m}R[T].$$

Hence there exist polynomials $h_i$ such that

$$P(T)^q = Q(T^q) - \sum_{i=1}^{\mu} x_i h_i(T). \tag{3}$$

From (2) it also follows that

$$P^q = \left(\sum_\nu r_{0,\nu} T^\nu\right)^q \tag{4}$$

by our choice of $q$. On the other hand

$$\bar{Q}(v^q) = \sum_\nu \bar{r}_{0,\nu}^q v^{q\nu} = \left(\sum_\nu \bar{r}_{0,\nu} v^\nu\right)^q$$

in $\kappa$. By definition of $\omega$ and (4) we then obtain that

$$\omega(\bar{Q}(v^q)) = P(c)^q. \tag{5}$$

Putting (5) together with (3) we obtain the desired equality. ∎

**3.6. Definition.** Let $(R,\mathfrak{m})$ be an Artinian local ring of mixed characteristic of length $l$, with algebraically closed residue field $\kappa$. Let notation be as in **(2.1)** and **(3.1)** and assume that $R$ is moreover Gorenstein. We define the map $\nabla\colon \kappa^l \to R$, which assigns to some tuple $u = (u_\alpha)_{\alpha\in\Delta_R}$, the element

$$\sum_{\alpha\in\Delta_R} \omega(u_\alpha)x^\alpha.$$

It follows from **(3.4)** that $\nabla$ is a bijection. (Caution: $\nabla$ is not additive and we merely have that $\nabla(k \cdot u) = \omega(k)\nabla u$, for $k \in \kappa$). We use the same notation $\nabla$ to indicate the map $(\kappa^l)^n \to R^n$ given by sending $(u_1,\ldots,u_n)$ to $(\nabla(u_1),\ldots,\nabla(u_n))$, where each $u_i \in \kappa^l$.

**3.7. Proposition.** *Let $(R,\mathfrak{m})$ be an Artinian local Gorenstein ring of mixed characteristic, with algebraically closed residue field $\kappa$. Let notation be as in **(2.1)** and **(3.1)**. Let $P \in R[T]$, where $T = (T_1,\ldots,T_n)$ are some variables. Then there exist $p_\alpha(U,V) \in \kappa[U,V]$, for $\alpha \in \Delta_R$, where $U = (U_{i\alpha} \mid i=1,\ldots,n, \alpha \in \Delta_R)$ and $V = (V_1,\ldots,V_m)$ are variables, and a quantifier free formula $\varphi(\bar{\boldsymbol{u}},\bar{\boldsymbol{v}})$ with parameters from $\kappa$, such that, for each Artinian local ring $(S,\mathfrak{n})$ of length $l$, containing $R$ and having an algebraically closed residue field $\lambda$ and for each $u = (u_{i\alpha}) \in \lambda^{nl}$, we can find $v \in \lambda^m$, such that $\varphi(u,v)$ holds in $\lambda$ and*

$$P(\nabla(u)) = \sum_{\alpha\in\Delta_R} \omega(p_\alpha(u,v))x^\alpha. \tag{1}$$

*Proof.* In order to prove this, we will use an induction on $\Delta_R$ to prove the following variant, for any $\alpha \in \Delta_R$.

**3.7.α. Claim.** *There exist $p_\delta^{(\alpha)} \in \kappa[U,V]$, for $\delta \in \Delta_R$ and $\delta \leq \alpha$, and $Q_\delta^{(\alpha)} \in R[U,V]$, for $\delta \in \Delta_R$ and $\alpha < \delta$, and there exists a quantifier free formula $\varphi^{(\alpha)}(\bar{\boldsymbol{u}},\bar{\boldsymbol{v}})$ with parameters from $\kappa$, such that, for each $S$ as in the statement of the Proposition and each $u \in \lambda^{ln}$, we can find an $v \in \lambda^m$, such that $\varphi^{(\alpha)}(u,v)$ holds in $\lambda$ and*

$$P(\nabla(u)) = \sum_{\delta\leq\alpha} \omega(p_\delta^{(\alpha)}(u,v))x^\delta + \sum_{\delta>\alpha} Q_\delta^{(\alpha)}(\omega(u,v))x^\delta. \tag{$1_\alpha$}$$

Hence Claim **(3.7.α)**, for $\alpha$ the largest element of $\Delta_R$ is precisely our original assertion.

*Proof of Claim.* The proof for $\alpha = 0$ is identical to the general case, so we will not treat it separately. Hence assume that **(3.7.α)** holds and let $\beta$ be the element of $\Delta_R$ which comes immediately after $\alpha$. (To treat the case $\beta = 0$ we do not have any quantifier free formula $\varphi^{(\alpha)}$ yet nor any $p_\delta^{(\alpha)}$). Take any $u \in \lambda^{ln}$ and let $v \in \lambda^m$ be the element whose existence is claimed by **(3.7.α)**. Let $w \in \lambda^{ln+m}$, such that (in our abbreviated notation) $w^q = (u, v)$. Let $\varphi^{(\beta)}(\bar{u}, \bar{v}, \bar{w})$ be the conjunction of $\varphi^{(\alpha)}(\bar{u}, \bar{v})$ together with the formula $\bar{w}^q = (\bar{u}, \bar{v})$. Apply Lemma **(3.5)** to the polynomial $Q_\beta^{(\alpha)}(U, V)$, to find $h_i \in R[U, V]$, for $i = 1, \ldots, \mu$, such that

$$(2) \qquad Q_\beta^{(\alpha)}(\omega(u, v)) = \omega(\bar{Q}_\beta^{(\alpha)}(u, v)) + \sum_{i=1}^{\mu} x_i h_i(\omega(w)).$$

Hence put $p_\delta^{(\beta)}$ equal to $p_\delta^{(\alpha)}$ in case $\delta \leq \alpha$ and to $\bar{Q}_\beta^{(\alpha)}$, for $\delta = \beta$. It is now clear that using (2) and $(1_\alpha)$, we can find $Q_\delta^{(\beta)} \in R[U, V, W]$ to make $(1_\beta)$ hold. ∎

**3.8. Theorem.** *The existentially closed models of the theory $\mathcal{A}\mathfrak{rt}_l$ are precisely the Artinian local Gorenstein rings of length $l$ with algebraically closed residue field.*

*Proof.* By **(1.6)** and **(2.5)**, we only need to show that an Artinian local Gorenstein ring $R$ of mixed characteristic and of length $l$ with algebraically closed residue field $\kappa$ is an existentially closed model of $\mathcal{A}\mathfrak{rt}_l$. Therefore, let $(S, \mathfrak{n}) \models \mathcal{A}\mathfrak{rt}_l$ be an overring of $R$, i.e., $R \subset S$. By enlarging $S$ if necessary, we may already assume that $S$ has algebraically closed residue field. We keep the notations of **(2.1)** and **(3.1)** for the ring $R$. Let $\varphi$ be an existential formula with parameters from $R$, which holds in $S$. We need to show that then already $R \models \varphi$. Since $\varphi$ is a conjunction of sentences of the form $\exists \bar{t} \psi(\bar{t})$, where

$$(1) \qquad \psi(\bar{t}) = \Big( \bigwedge_{j < M} P_j(\bar{t}) = 0 \Big) \wedge \Big( \bigwedge_{M \leq j < N} P_j(\bar{t}) \neq 0 \Big),$$

with $P_j(T) \in R[T]$ and $T = (T_1, \ldots, T_n)$, we may without loss of assumption assume that $\varphi$ is already of the form $\exists \bar{t} \psi(\bar{t})$. Apply **(3.7)** to each polynomial $P_j$, for $j < N$, to find quantifier free formulae $\varphi_j(\bar{u}, \bar{v})$ with parameters from $\kappa$, and polynomials $p_{j\alpha}(U, V) \in \kappa[U, V]$, for $\alpha \in \Delta_R$, where $U = (U_{i\alpha} \mid i = 1, \ldots, n, \alpha \in \Delta_R)$ and $V = (V_1, \ldots, V_m)$ are variables, so that for each $u = (u_{i\alpha}) \in \lambda^{nl}$, we can find an $v \in \lambda^m$, such that $\varphi_j(u, v)$ holds in $\lambda$ and

$$(2) \qquad P_j(\nabla(u)) = \sum_{\alpha \in \Delta_R} \omega(p_{j\alpha}(u, v)) x^\alpha.$$

Now, let $t = (t_1, \ldots, t_n) \in S^n$, such that $S \models \psi(t)$. Using that $\mathcal{E}_R$ is a Witt basis for $S$ according to **(3.4)**, we can find $u \in \lambda^{ln}$, such that $\nabla(u) = t$. Let $v \in \lambda^m$ be such that (2) becomes true. Let $\xi(\bar{u}, \bar{v})$ be the formula

$$\Big( \bigwedge_{j < N} \varphi_j \Big) \wedge \Big( \bigwedge_{j < M} \bigwedge_{\alpha \in \Delta_R} p_{j\alpha}(\bar{u}, \bar{v}) = 0 \Big) \wedge \Big( \bigwedge_{M \leq j < N} \bigvee_{\alpha \in \Delta_R} p_{j\alpha}(\bar{u}, \bar{v}) \neq 0 \Big).$$

Using once more that $\mathcal{E}_R$ is a Witt base of $S$, we see that $\lambda \models \xi(u, v)$. Since $\kappa$ is algebraically closed, we can find $u_0 \in \kappa^{ln}$ and $v_0 \in \kappa^m$, such that $\kappa \models \xi(u_0, v_0)$, which means that $R \models \psi(t_0)$, where $t_0 = \nabla(u_0)$, as required. ∎

## 4. Model Completeness

**4.1. Lemma.** *The class of all equicharacteristic local rings of depth zero is elementary.*

*Proof.* Consider, for each $n \in \mathbb{N}_0$, the following sentence

$$(4.1.1) \qquad \mathtt{Ec}_n \stackrel{\mathrm{def}}{=} \mathtt{Loc} \wedge (\exists \boldsymbol{x})[\mathtt{Mass}(\boldsymbol{x}) \wedge (n \cdot \mathbf{1} = 0 \leftrightarrow n \cdot \boldsymbol{x} = \mathbf{0})].$$

If $\mathtt{Ec}_n$ holds in a ring $R$, then this ring is local, of depth zero and $n$ is zero in $R$, if and only if, it is zero in the residue field of $R$, since the annihilator of $x$ is precisely the maximal ideal. So adding all $\mathtt{Ec}_n$ to $\mathcal{R}\mathrm{ing}$ gives the theory of all equicharacteristic local rings of depth zero. ∎

**4.2. Lemma.** *The class of all local rings of depth zero with algebraically closed residue field is elementary.*

*Proof.* It is well-known that the class of algebraically closed fields is elementary: except for the field axioms, one needs to add, for each degree $d \geq 1$, a sentence expressing that each monic polynomial of degree $d$ over this field has a root within the field. Using the techniques of **(1.1)** we can easily deduce the required sentences from this by setting

$$(4.2.1) \qquad \mathtt{Root}_d(\boldsymbol{x}) \stackrel{\text{def}}{=} (\forall \bar{\boldsymbol{a}})(\exists \boldsymbol{y})[\boldsymbol{x}(\boldsymbol{y}^d + \sum_{i=1}^d \boldsymbol{a_i}\boldsymbol{y}^{d-i}) = \boldsymbol{0}],$$

where we have replaced 'being zero' by 'being annihilated by $x$', the latter being equivalent under the assumption $\mathtt{Mass}(x)$ with 'lying in the maximal ideal'. Therefore, adding

$$(4.2.2) \qquad \{\mathtt{Loc}\} \cup \{ (\exists \boldsymbol{x})[\mathtt{Mass}(\boldsymbol{x}) \wedge \mathtt{Root}_d(\boldsymbol{x})] \mid d \geq 1 \}$$

to $\mathcal{R}\mathfrak{ing}$ gives the desired theory. ∎

**4.3. Theorem.** *The theory $\mathcal{G}\mathfrak{or}_l$ of all Artinian local Gorenstein rings of length $l$ having an algebraically closed residue field is model-complete.*

*Proof.* Let

$$(4.3.1) \qquad \mathtt{Min}(\boldsymbol{x}) \stackrel{\text{def}}{=} (\forall \boldsymbol{a})(\exists \boldsymbol{b})[\boldsymbol{x} \neq \boldsymbol{0} \wedge (\boldsymbol{a} = \boldsymbol{0} \vee \boldsymbol{x} = \boldsymbol{ba})].$$

For any model $R$ of $\mathcal{A}\mathfrak{rt}_l$, we have that $(\exists \boldsymbol{x})\mathtt{Min}(\boldsymbol{x})$ holds in $R$, if and only if, there is a unique minimal ideal in $R$, which is equivalent with saying that $R$ is Gorenstein. Together with **(4.2)**, we then obtain the desired first order theory $\mathcal{G}\mathfrak{or}_l$. Let $\mathbb{K}_l$ denote the class of all models of $\mathcal{G}\mathfrak{or}_l$. As $\mathbb{K}_l$ is a subclass of the class of all Artinian local rings of length $l$, any member of $\mathbb{K}_l$ which is existentially closed in the latter class is already so in $\mathbb{K}_l$. In other words, by **(3.8)**, each member of $\mathbb{K}_l$ is existentially closed. By [**Hod 2**, Theorem 8.3.1] this implies that its theory $\mathcal{G}\mathfrak{or}_l$ is model-complete. ∎

*Remark.* Using **(4.1)**, one shows that also the class of $l$-dimensional local Gorenstein algebras over an algebraically closed field is model-complete.

**4.4. Theorem.** *The theory $\mathcal{A}\mathfrak{rt}_l$ of Artinian local rings of length at most $l$ is companionable.*

*Proof.* This follows immediately from [**Hod 2**, Theorem 8.3.6] and **(4.3)**. ∎

## 5. A Counterexample to Quantifier Elimination

The following example shows that Quantifier Elimination does not hold for $\mathcal{G}\mathfrak{or}_l$, even for $l = 2$. Note that $l = 2$ is actually quite more favourable than the other cases: if we add to $\mathcal{G}\mathfrak{or}_2$ an axiom fixing the characteristic to be either zero or prime (implying the equicharacteristic case), then we obtain even a complete theory. Indeed, any model is then of the form $\kappa[T]/(T^2)$, with $\kappa$ an algebraically closed field of the prescribed characteristic. This is still true for $l = 3$, but not for $l > 3$. Nonetheless even in this simplest case Quantifier Elimination already fails.

To exhibit a counterexample, we will show that amalgamation over substructures fails for $\mathcal{G}\mathfrak{or}_2$[5] (see [**Hod 2**, Theorem 8.4.1]). To this end, we need to take a closer look at automorphisms. We will not provide proofs (they are not that hard), for they will appear in a subsequent paper [**Sch 2**]. The easiest case is that of $\kappa$ an algebraically closed field of prime characteristic.[6] Then any automorphism of $\kappa[T]/(T^2)$, where $T$ is one variable, is of the form $\Lambda(\sigma, d)$, where $\sigma \in \mathrm{Aut}(\kappa)$ and $d \in \kappa^\times$, defined by

$$\Lambda(\sigma, d)\colon u + vT \mapsto \sigma(u) + d\sigma(v)T.$$

In fact, any embedding of $\kappa[T]/(T^2)$ into $\lambda[T]/(T^2)$, where $\lambda$ is an extension of $\kappa$, is of the form $\Lambda(\sigma, d)$, with $\sigma \in \mathrm{Aut}(\kappa)$ and $d \in \lambda^\times$.

---

[5]If $l > 3$ then even simpler counterexamples to Quantifier Elimination for $\mathcal{G}\mathfrak{or}_l$ exist. For instance, one cannot amalgamate $\kappa[X]/(X^4)$ and $\kappa[X, Y]/(X^2, Y^2)$ over $\kappa$.

[6]If the characteristic is zero, different automorphisms might exist, due to the existence of non-zero derivations, obstructing even further any possible amalgamation.

**Counterexample.** Let $\omega$ be the finite field $\mathbb{F}_p$ with $p$ elements (where $p$ is a prime). Let $\kappa$ be the algebraic closure of $\omega$. Let $V = \omega[X,Y]/(X,Y)^2$ and let $R = \kappa[T]/(T^2)$, so that $R$ is a model of $\mathcal{G}\mathfrak{or}_2$.

Consider the following two embeddings $\phi_k$ of $V \hookrightarrow R$, where $k = 1, 2$. Let $\phi_k$ send an element $r = u + vX + wY$ to $u + (v\theta_k + w)T$, where $u, v, w \in \omega$ and $\theta_k \in \kappa \setminus \omega$. (Note that these are indeed embeddings since $\theta_k \notin \omega$.) We claim no amalgam in $(\mathcal{G}\mathfrak{or}_2)_\forall$ over these two embeddings can exist in general. Suppose the contrary and assume there exists an amalgam $S$. Without loss of generality, we may assume that $S \models \mathcal{G}\mathfrak{or}_2$, so that $S = \lambda[T]/(T^2)$, with $\lambda$ algebraically closed. Hence there exist $\sigma_k \in \text{Aut}(\kappa)$ and $d_k \in \lambda^\times$, such that the following diagram commutes

$$\begin{array}{ccc} V & \xrightarrow{\phi_1} & R \\ \phi_2 \downarrow & & \downarrow \Lambda(\sigma_1, d_1) \\ R & \xrightarrow{\Lambda(\sigma_2, d_2)} & S. \end{array}$$

This means that, for all $v, w \in \omega$, we must have that

$$d_1(v\sigma_1(\theta_1) + w) = d_2(v\sigma_2(\theta_2) + w).$$

In other words, $d_1 = d_2$ and $\sigma_1(\theta_1) = \sigma_2(\theta_2)$. In particular, we see that $\theta_1$ and $\theta_2$ must be conjugate (i.e., have the same minimal polynomial), which obviously does not necessarily hold.

FIELDS INSTITUTE
222 COLLEGE STREET
TORONTO, ONTARIO, M5T 3J1 (CANADA)
*E-mail address*: hschoute@fields.utoronto.ca